\documentclass[a4paper,10pt]{amsart}
\usepackage{amsmath}
\usepackage{amsthm}
\usepackage{amssymb}
\usepackage{amsfonts}
\usepackage{amsxtra}
\usepackage[mathscr]{eucal}
\DeclareMathOperator{\ord}{ord}
\theoremstyle{definition}
\newtheorem{thm}{Theorem}[section]
\newtheorem{lem}[thm]{Lemma}
\newtheorem{cor}[thm]{Corollary}

\newtheorem{defi}[thm]{Definition}

\begin{document}
\title{Primitive Divisors of some Lehmer-Pierce Sequences}
\author{Anthony Flatters}
\address{School of Mathematics, University of East Anglia, Norwich NR4 7TJ, UK}
\email{Anthony.Flatters@uea.ac.uk}
\begin{abstract}
We study the primitive divisors of the terms of $(\Delta_n)_{n\geqslant 1}$, where $\Delta_n=N_{K/\mathbb{Q}}(u^n-1)$ for $K$ a real quadratic field, and $u>1$ a unit element of its ring of integers. The methods used allow us to find the terms of the sequence that do not have a primitive prime divisor.
\end{abstract}
\maketitle

\section{Introduction}
Let $A=(a_n)_{n\geqslant 1}$ be an integer sequence. A prime $p$ dividing a term $a_n$ is called a {\em primitive prime divisor} (PPD for short) of $a_n$ if $p$ does not divide $a_m$ for any $m<n$ with $a_m\neq 0$. Sequences whose terms all have primitive divisors beyond some point are of great interest in number theory. 
\begin{defi}
Let $A=(a_n)_{n\geqslant 1}$ be an integer sequence. Define 
\begin{displaymath}
Z(A)=\max \{n:a_n \textrm{ does not have a primitive prime divisor}\}
\end{displaymath}
if this set is finite, otherwise set $Z(A)=\infty$. The number $Z(A)$ is called the \emph{Zsigmondy Bound} for the sequence $A$. \end{defi}
In \cite{Bang}, Bang considered the sequence $(a^n-1)_{n\geqslant 1}$, where $1<a\in\mathbb{Z}$ and showed that $Z((a^n-1)_{n\geqslant 1})\leqslant 6$. Zsigmondy in \cite{Zsigmondy} proved the more general result that given any positive coprime integers $a,b$ with $a>b$, the sequence $(a^n-b^n)_{n\geqslant 1}$ has a primitive prime divisor for all terms beyond the sixth. The sequence studied by Zsigmondy satisfies a binary linear recurrence relation, and much of the work in this area has concentrated on these types of sequences. In \cite{Carmichael}, Carmichael showed that for any real Lucas or Lehmer sequence $L$, $Z(L)\leqslant 12$. Carmichael's result was later completed by Bilu, Hanrot and Voutier, and in \cite{Bilu} they showed, using powerful methods from transcendence theory and computational number theory, that for any Lucas or Lehmer sequence $L^{\prime}$, $Z(L^{\prime})\leqslant 30$. Moreover, they were able to explicitly describe all Lucas and Lehmer numbers without a primitive divisor and hence show that this bound is sharp.\newline\newline Many arithmetic properties of linear recurrence sequences have analogues for elliptic recurrence sequences. In \cite{Silverman}, it is shown that if $E$ is an elliptic curve in Weierstrass form defined over $\mathbb{Q}$, and $P\in E(\mathbb{Q})$ is a non-torsion point, then the associated elliptic divisibility sequence (the denominators of the $x$-coordinates of $nP$) has a finite Zsigmondy bound. For elliptic curves in global minimal form, it seems likely that this bound is uniform, and the papers \cite{EverestWard4}, \cite{Ingram} exhibit infinite families of elliptic curves with a uniform Zsigmondy bound. \newline\newline The result of Zsigmondy can be generalised to a number field setting, where $a,b$ are now algebraic integers of a number field $K$, so $a^n-b^n$ lies in the ring of integers $R$, of $K$. The principal ideal $(a^n-b^n)$ has a factorisation into a product of prime ideals of $R$, which is unique. Therefore, we can ask which terms of a sequence $S$ of algebraic integers have a primitive prime ideal divisor (or PPID for short), i.e. for which $n$ is there a prime ideal $\mathfrak{p}$ which divides the $n$th term, but not any preceding term. We therefore define the Zsigmondy bound $Z_I(S)$, to be the maximal value of $n$ for which the $n$th term of the sequence does not have a PPID. 
\newline\newline In Schinzel's paper \cite{Schinzel}, he proved the following theorem;
\begin{thm}[Schinzel] \label{Schinzelthm}
Let $A,B$ be coprime integers of an algebraic number field such that $\frac{A}{B}$ is not a root of unity. Then the expression $A^n-B^n$ has a PPID for all $n>n_0(d)$, where $d$ is the degree of the extension $\mathbb{Q}\left ( \frac{A}{B} \right )/\mathbb{Q}$. 
\end{thm}
So, for these sequences the Zsigmondy bound $Z_I$ is finite and an easy corollary of Schinzel's theorem is the following.
\begin{cor}\label{Ucorollary}
Let $K$ be a real quadratic field, $R$ its ring of integers, and let $\alpha\in R\setminus\{\pm 1\}$ be a unit. Let $f$ denote the minimum polynomial of $\alpha$ over $\mathbb{Q}$ and define the integer sequence, $\Delta=(\Delta_n(f))_{n\geqslant 1}$, by setting \begin{displaymath}
\Delta_n(f)=N_{K/\mathbb{Q}}(\alpha ^n -1)\textrm{.}  \end{displaymath}
 Then there exists a positive integer $C_1$, so that for all units $\alpha$ of norm 1, $Z(\Delta)\leqslant C_1$. There exists a positive integer $C_2$ such that for all units $\alpha$ of norm $-1$, $\Delta_n(f)$ has a primitive prime divisor for any $n>C_2$ with $n\not\equiv 2 \pmod{4}$. \end{cor}
The sequence $\Delta$, for a general algebraic integer $\alpha$, was examined by Pierce in his paper \cite{Pierce}, where he looked at what form the factors of $\Delta_n(f)$ take and what conditions are necessary for the congruence $f(x)\equiv 0 \pmod{p}$, where $p$ is a prime, to have a solution $x\in\mathbb{F}_p$. In \cite{Lehmer}, Lehmer developed a deeper insight into the factors of the terms $\Delta_n(f)$, and applied this information to show that certain $\Delta_n(f)$ were prime. Lehmer was interested in the growth rate of the sequence $\Delta$, and he remarked that if none of the roots of $f$ had absolute value 1, then $\frac{\Delta_n(f)}{\Delta_{n-1}(f)}$ converges, and $M(f)$ was written for the limit. For his purposes, polynomials with small values of $M(f)$ were desirable; in \cite{EverestWard2} a heuristic argument is put forward that suggests the density of primes in $\Delta$ is proportional to $\frac{1}{M(f)}$. We therefore say that for $\alpha$ an algebraic integer, the sequence $\Delta$, defined in Corollary \ref{Ucorollary}, is called the {\em Lehmer-Pierce sequence} associated to $\alpha$. The sequence $\Delta$ is also of interest in algebraic dynamics, since to $f$ there is an associated matrix called the companion matrix and multiplication by this matrix induces an endomorphism $E:\mathbb{T}^N\longrightarrow \mathbb{T}^N$. When none of the roots of $f$ have absolute value 1, $E$ is an {\em ergodic} transformation with respect to Lebesgue measure, $\vert \Delta_n(f) \vert $ counts the number of points of period $n$ under $E$, and the topological entropy of $E$ is equal to $\log{M(f)}$. A much more detailed account of the connection between $M(f)$ and dynamical systems can be found in \cite{EverestWard3}.\newline\newline The sequence $\Delta$ also has some combinatorial applications. For example, when $u=1+\sqrt{2}$, $\vert \Delta_n \vert $ counts the number of $2 \times 2$ tiles in all tilings of a $3 \times (n+1)$ rectangle with $1 \times 1$ and $2 \times 2$ square tiles; more details about this sequence are provided on Sloane's website \cite[A095977]{EISW}. Similarly, when $u=\frac{3+\sqrt{5}}{2}$, $\Delta_n$ appears in combinatorics - see \cite[A004146]{EISW}. In addition, certain quadratic Lehmer-Pierce sequences count the sizes of groups: the groups being $E(\mathbb{F}_{p^n})$, where $E$ is a given elliptic curve and $p$ is a fixed prime. \newline\newline In this article our aim is to find the numbers $C_1,C_2$ from Corollary \ref{Ucorollary} associated to the sequence $\Delta=(\Delta_n)_{n\geqslant 1}=(N_{K/\mathbb{Q}}(u^n-1))_{n\geqslant 1}$, where $K$ is a real quadratic field and $u$ is a fixed unit in its ring of integers. 
\begin{thm}
Let $K$ be a quadratic field, $\alpha\neq \pm 1$ a positive quadratic unit, and let $\Delta$ be the Lehmer-Pierce sequence associated to $\alpha$. Then for each $\alpha$ of norm 1, $\Delta$ has a primitive prime divisor for all terms beyond the twelfth. For each $\alpha$ of norm $-1$, then for $n>24$, $\Delta_n$ fails to have a primitive prime divisor if and only if $n\equiv 2 \pmod{4}$.
\end{thm} 
It is easy to see that when $u$ has norm 1, $\Delta$ satisfies a ternary linear recurrence relation, and when $u$ has norm $-1$, a quaternary linear recurrence - see \cite{EverestWard5}. In addition, it is remarked that it seems likely that when $u=2+\sqrt{3}$, $Z(\Delta)=6$, and in our later discussion we verify that this is indeed the case. To date, not much is known about primitive prime divisors of the terms $\Delta_n$ for arbitrary algebraic integers $\alpha$, and it would be interesting to know which other Lehmer-Pierce sequences have the property that $Z(\Delta)$ is finite.

\section{A Criterion for Primitive Divisor Failure} 
We begin with a proof of Corollary \ref{Ucorollary} as it will be instrumental in obtaining a condition that will need to be satisfied if $\Delta_n$ fails to have a PPD.

\begin{proof}[Proof of Corollary~\ref{Ucorollary}]
Define $A_n=\alpha ^n-1$ and $B_n=\beta ^n-1$, where $\beta$ is the algebraic conjugate of $\alpha$. There are only two ways in which $\Delta_n$ could fail to have a primitive prime divisor, and they are the following:
\begin{enumerate}
\item Both $A_n$ and $B_n$ fail to have PPIDs;
\item Every PPID of $A_n$ has already appeared before as a divisor of $B_m$ for some $m<n$.
\end{enumerate}
Suppose then that $\mathfrak{p}$ is a PPID of $A_n$ but that $\mathfrak{p}\vert B_m$ for some $m<n$. Then \begin{displaymath}
(\beta ^m-1)=\mathfrak{p}\mathfrak{q} \end{displaymath}
for some integral ideal $\mathfrak{q}$. Hence, multiplying through by $(\alpha ^m)$, \begin{displaymath}
(\alpha ^m)(\beta ^m-1)=\mathfrak{p}\mathfrak{q}\textrm{.} \end{displaymath}
If $\alpha$ has norm 1, this therefore implies that $\mathfrak{p}$ divides $A_m$, which cannot be the case as $\mathfrak{p}$ is a PPID of $A_n$. \newline If $\alpha$ has norm $-1$ and $m$ is even, then by the same method as above we can deduce that possibility 2 will not occur. If $\alpha$ has norm $-1$ and $m$ is odd, a slightly different argument is needed. If possibility 2 occurs in this case, we have that $\mathfrak{p}\vert (\alpha ^m+1)$. Therefore, $\alpha^m\equiv -1  \pmod{\mathfrak{p}}$ and so $\alpha ^{2m}\equiv 1 \pmod {\mathfrak{p}}$. Now as $\mathfrak{p}$ is a primitive divisor of $A_n$, $\alpha$ has order $n$ in the group $ \left ( R/\mathfrak{p}\right )^{\ast}$. Therefore $n\vert 2m$. Since $m<n$, this is enough to secure that $n=2m$, and we conclude that possibility 2 can only hold in the case when $n$ is twice an odd integer. If $n\equiv 2 \pmod{4}$, then $n=2k$ for some odd integer $k$ and in this case $\Delta_n=-\Delta_k^2$, so $\Delta_n$ can never have any primitive prime divisors. We have deduced that if $\Delta_n$ fails to have a PPD, then both $A_n$ and $B_n$ fail to have PPIDs except in the case where $\alpha$ has norm $-1$ and then all terms which satisfy property 2, are those with $n\equiv 2 \pmod{4}$. The fact that $\Delta_n$ fails to have a primitive divisor beyond some point if $n\equiv 2 \pmod{4}$ was first pointed out by Gy\"ory.
\newline\newline Hence for units of norm 1, $\Delta_n$ will only fail to have a PPD, when condition 1 holds. So by Theorem \ref{Schinzelthm}, this tells us that $Z_I((A_n)_{n\geqslant 1})<c_1$, and $Z_I((B_n)_{n\geqslant 1})<c_2$, where $c_1,c_2$ are uniform constants, and so for all units $\alpha$ of norm 1, $Z(\Delta)$ is uniformly bounded. If $\alpha$ has norm $-1$, then $\Delta_n$ will fail to have a PPD when $n\equiv 2 \pmod{4}$ and when condition 1 holds. Applying Theorem \ref{Schinzelthm} again gives the required result.
\end{proof}
From now on, $K$ denotes a real quadratic field we will write $N$ for the field norm $N_{K/\mathbb{Q}}$.
\begin{lem}
Let $u\in R\setminus \{\pm 1 \}$ be a quadratic unit of norm 1. Then for any $n>6$, if $\Delta_n$ fails to have a primitive prime divisor we have
\begin{equation} \label{eq:eq1}
N(\phi_n(u))\Bigl \vert n^2\textrm{,} \end{equation}
where $\phi_n(x)\in \mathbb{Z}[x]$ denotes the $n$th cyclotomic polynomial. Moreover, if $u$ has norm $-1$ then for any $n>6$ with $n \not \equiv 2 \pmod{4}$, if $\Delta_n$ fails to have a PPD then (\ref{eq:eq1}) holds.
\end{lem}
\begin{proof}
 Apply Lemma 4 of \cite{Schinzel} to deduce that if $\mathfrak{p}$ is not a PPID of $A_n$ or $B_n$, then for $n>6$,  \begin{displaymath} \ord_{\mathfrak{p}}(\phi_n(u))\leqslant \ord_{\mathfrak{p}}(n) \end{displaymath} and \begin{displaymath} \ord_{\mathfrak{p}}(\phi_n(v))\leqslant \ord_{\mathfrak{p}}(n)\textrm{.} \end{displaymath} 
Adding these two inequalities together tells us that 
\begin{displaymath} 
\ord_{\mathfrak{p}}(N(\phi_n(u)))\leqslant \ord_{\mathfrak{p}}(n^2) \textrm{,} \end{displaymath}
and so we have proved the statement of the Lemma.
\end{proof}

Using (\ref{eq:eq1}), we can express this result in a way that will allow us to obtain an upper bound on $n$ such that $\Delta_n$ has no PPD.

\begin{thm} \label{asymptoticthm}
Let $1<u\in R$ be a unit, and $6<n\in\mathbb{N}$. If $u$ has norm 1 and $\Delta_n$ has no primitive prime divisor, then
\begin{equation} \label{eq:eq2}
\log{n} -2\log{\log{n}} -\frac{4}{\log{n}}<2.02819-\log{\log{u}}\textrm{.}
\end{equation}
If $u$ has norm $-1$, $n \not \equiv 2 \pmod{4}$, and $\Delta_n$ has no primitive prime divisor, then
\begin{equation} \label{eq:eq3}
\log{n} -2\log{\log{n}} -\frac{4}{\log{n}}<2.71072-\log{\log {u}} \textrm{.} \end{equation}
\end{thm}

\begin{proof}
Recall the factorisation of $x^n-1$ into a product of cyclotomic polynomials as follows
\begin{displaymath}
x^n-1=\prod_{d \vert n}\phi_d (x)\textrm{.} \end{displaymath}
Hence we have the following factorisation of $\Delta_n$
\begin{displaymath} \lvert \Delta_n \rvert =\prod_{d \vert n}\lvert N(\phi_d (u)) \rvert \textrm{.} \end{displaymath}
Taking logarithms now gives
\begin{displaymath}
\log(\lvert N(u^n-1)\rvert)=\sum_{d \vert n} \log(\lvert N(\phi_d(u))\rvert)\textrm{.} \end{displaymath}
Applying the M\"obius Inversion Formula for arithmetical functions now yields 
\begin{eqnarray}
\ \label{eq:eq4} \log(\lvert N(\phi_n(u))\rvert)= \sum_{d \vert n}\log(\lvert N(u^d-1)\rvert )\mu\left(\frac{n}{d}\right)\textrm{.}  \end{eqnarray}

Now using (\ref{eq:eq4}), we are going to estimate the size of $\lvert N(\phi_n(u)) \rvert$. If $u$ is a unit of norm 1, then 
\begin{eqnarray}
\log{\lvert N(u^d-1) \rvert}  &=& \log{\lvert u^d-1 \rvert} +\log{\lvert v^d-1 \rvert}  \nonumber \\
\ &=& \log{\lvert u^d-1 \rvert} + \log{\left \lvert \frac{1}{u^d} -1 \right \rvert} \nonumber \\
\ &=& \log{\lvert u^d \rvert} + 2\log \left \lvert 1-\frac{1}{u^d} \right \rvert \textrm{.} \nonumber
\end{eqnarray}
Therefore, by (\ref{eq:eq4}) we have
\begin{eqnarray}
\log(\lvert N(\phi_n(u))\rvert) &=& \sum_{d \vert n}\log{\lvert u^d \rvert}\mu\left (\frac{n}{d}\right ) + 2\sum_{d \vert n}\log \left \lvert 1-\frac{1}{u^d} \right \rvert \mu \left ( \frac{n}{d} \right ) \nonumber \\
\ &=& \phi(n)\log{u} +2\sum_{d \vert n}\log \left \lvert 1-\frac{1}{u^d} \right \rvert \mu \left ( \frac{n}{d} \right ) \textrm{.} \nonumber 
\end{eqnarray}
Define $S:=2\sum_{d \vert n}\log \left \lvert 1-\frac{1}{u^d} \right \rvert \mu \left ( \frac{n}{d} \right )$. Using the Taylor expansion for $\log(1-x)$, we obtain that 
\begin{displaymath}
\vert S \vert = 2 \left \lvert \sum_{m=1}^\infty\frac{1}{m}\sum_{d \vert n}\frac{1}{u^{md}} \mu \left ( \frac{n}{d} \right ) \right \rvert \textrm{.} \end{displaymath}
Hence,
\begin{eqnarray}
\ \vert S \vert &<&  2\sum_{m=1}^\infty\frac{1}{m}\sum_{d=1}^\infty\frac{1}{u^{md}} \nonumber \\
\ &=& 2\sum_{m=1}^\infty \frac{1}{m}\left ( \frac{u^{-m}}{1-u^{-m}}\right ) \textrm{.} \nonumber
\end{eqnarray}
Since $u$ has norm 1, $u\geqslant \frac{3+\sqrt{5}}{2}$. In addition, $m\geqslant 1$ so 
\begin{equation}\label{eq:eq5}
\vert S \vert <3.23607\sum_{m=1}^\infty \frac{1}{m\left (\frac{3+\sqrt{5}}{2} \right )^m}\textrm{.} \end{equation}
The sum in (\ref{eq:eq5}) is equal to $-\log\left ( 1-\frac{2}{3+\sqrt{5}} \right )$, and so
\begin{displaymath}
\vert S \vert <1.55724 \textrm{,} \end{displaymath}
which therefore yields that
\begin{equation}
 \log(\lvert N(\phi_n(u))\rvert)>\phi(n)\log{u}-1.55724 \textrm{.} \nonumber \end{equation}

Now we use the fact that if $\Delta_n$ has no PPDs, then $\lvert N(\phi_n(u))\rvert\leqslant n^2$.
Therefore, we have the following relation
\begin{equation}
 \label{eq:eq7} u^{\phi(n)}<e^{1.55724}n^2\textrm{.}  \end{equation}
Taking logarithms twice of both sides we obtain
\begin{equation}
 \log(\phi(n))+\log{\log{u}}< \log(1.55724+2\log{n}) \textrm{.} \nonumber
\end{equation}
Since $n>6$, $\log{n}>1$, hence we have that
\begin{displaymath}
\log{n}+\sum_{p \vert n}\log \left (1-\frac{1}{p} \right )< \log(3.55724)-\log{\log{u}} +\log{\log{n}} \textrm{,} \end{displaymath}
and therefore
\begin{displaymath}
\log{n} < 1.26899-\log{\log{u}}+\log{\log{n}}-\sum_{p\vert n}\log\left (1-\frac{1}{p}\right )\textrm{.} \end{displaymath}
Noting now that for all primes $p$, $-\log\left ( 1-\frac{1}{p} \right )\leqslant \frac{1}{p}+\frac{1}{p^2}$ yields that \begin{displaymath}
\log{n} < 1.26899 -\log{\log {u}} +\log{\log{n}}+ \sum_{p\vert n} \frac{1}{p} + \sum_{p\vert n}\frac{1}{p^2} \textrm{.} \end{displaymath}
By Proposition 2.3.3, page 72 in \cite{Jameson}, the last term in our previous inequality is at most $\log(\zeta(2))$, where $\zeta(s)$ denotes the Riemann-Zeta function. Therefore,
\begin{displaymath}
\log{n}< 1.76669-\log{\log {u}} +\log{\log{n}}+\sum_{p\leqslant n}\frac{1}{p}\textrm{.}\end{displaymath}
In \cite{Villarino}, the following estimate is derived \begin{displaymath}
\sum_{p\leqslant n}\frac{1}{p} < \log{\log{n}}+B+\frac{4}{\log{n}}\textrm{,} \end{displaymath}
where $B$ is a numerical constant whose value is approximately equal to 0.2614972128.
Inserting all this information into our inequality yields
\begin{displaymath}
\log{n}-2\log{\log{n}} - \frac{4}{\log{n}} < 2.02819 - \log{\log {u}} \textrm{.}  \end{displaymath}
If $u$ is a unit of norm $-1$, then \begin{displaymath}
\log\vert N(u^d-1) \vert = \log \vert u^d \vert +\log \left \lvert 1-\frac{1}{u^d} \right \rvert + \log\left \lvert 1-\frac{(-1)^d}{u^d} \right \rvert \textrm{.} \end{displaymath}
Plugging this in to equation (\ref{eq:eq4}), we have
\begin{eqnarray}
\log(\lvert N(\phi_n(u))\rvert)&=&\sum_{d \vert n} \log \vert u^d \vert \mu \left ( \frac{n}{d} \right ) +\sum_{d\vert n} \log \left \lvert 1-\frac{1}{u^d} \right \rvert \mu \left ( \frac{n}{d} \right )+\sum_{d\vert n} \log \left \lvert 1-\frac{(-1)^d}{u^d} \right \rvert\mu \left ( \frac{n}{d} \right ) \nonumber \\
\ &=& \phi(n)\log{u} + \sum_{d\vert n} \log \left \lvert 1-\frac{1}{u^d} \right \rvert \mu \left ( \frac{n}{d} \right )+\sum_{d\vert n} \log \left \lvert 1-\frac{(-1)^d}{u^d} \right \rvert\mu \left ( \frac{n}{d} \right )\textrm{.} \nonumber
\end{eqnarray}
Define $S_1=\sum_{d\vert n} \log \left \lvert 1-\frac{1}{u^d} \right \rvert \mu \left ( \frac{n}{d} \right )$ and $S_2=\sum_{d\vert n} \log \left \lvert 1-\frac{(-1)^d}{u^d} \right \rvert\mu \left ( \frac{n}{d} \right )$. Again, using the Taylor expansion for $\log(1-x)$, and estimating these sums in the same way we did for $S$, we get
\begin{eqnarray}
\vert S_i \vert & \leqslant &  \sum_{m=1}^\infty \frac{1}{m} \sum_{d\vert n } \frac{1}{u^{md}} \nonumber \\
\ &<& \sum_{m=1}^\infty \frac{1}{m} \left ( \frac{u^{-m}}{1-u^{-m}}\right ) \textrm{.} \nonumber 
\end{eqnarray}
Noting that since $u$ is a unit of norm $-1$, $u\geqslant \frac{1+\sqrt{5}}{2}$, we see that
\begin{displaymath}
\vert S_i \vert <2.61804\sum_{m=1}^\infty \frac{1}{m\left (\frac{1+\sqrt{5}}{2} \right )^m}\textrm{.} \end{displaymath}
Once again, this sum is equal to $-\log\left ( 1-\frac{2}{1+\sqrt{5}} \right )$, thus
\begin{displaymath}
\vert S_i \vert <2.51966 \textrm{,} \end{displaymath}
and so 
\begin{displaymath} 
\log(\lvert N(\phi_n(u))\rvert)>\phi(n)\log{u}-5.03933 \textrm{.} \end{displaymath}
 Exponentiating this relation, we arrive at 
\begin{equation} \label{eq:eq8} \lvert N(\phi_n(u)) \rvert >\frac{u^{\phi(n)}}{e^{5.03933}}\textrm{.} \end{equation}
Running through the same calculation as before gives us the desired inequality.
\end{proof}

\section{Main Results}

\subsection{Units of Norm 1}
If $u>1$ is a unit of norm 1, then $u\geqslant (\frac{1+\sqrt{5}}{2})^2$. Inserting this into (\ref{eq:eq2}), we have that if $\Delta_n$ has no PPD, then \begin{displaymath}
\log{n}-2\log{\log{n}} - \frac{4}{\log{n}} <2.06650 \end{displaymath}

It is now clear that $n$ is bounded, since $g(x):=\log{x}-2\log{\log{x}}-\frac{4}{\log{x}}$ is an increasing function on $(e,\infty)$. Then, since $g(n)$ is bounded above, $n$ is also and so using Maple 9.5 to solve $g(x)=2.06650$ we conclude that 
\begin{displaymath}
n \leqslant 604 \textrm{.} \end{displaymath}
We can now improve this further because we know that inequality (\ref{eq:eq7}) must be satisfied in order that $\Delta_n$ has no PPD. We also know that $u\geqslant \left( \frac{1+\sqrt{5}}{2} \right )^2$. So we do a case by case check of the values of $n$ between 7 and 604 for which 
\begin{equation}
\label{eq:eq9} \left( \frac{1+\sqrt{5}}{2} \right )^{2\phi(n)}-e^{1.55724}n^2<0\textrm{.} \end{equation}
Instructing Maple 9.5 to compute the left hand side of the above inequality for each $n$ in our range and observing when the quantity is negative yields that
\begin{displaymath} n\leqslant 30\textrm{.} \end{displaymath}
More precisely, inequality (\ref{eq:eq9}) only holds when $n=8,9,10,12,14,18,24$ or 30. Now a bare hands approach is required to see if we can lower the bound.\newline\newline If we choose $u$ so that $u\geqslant C>\left( \frac{1+\sqrt{5}}{2} \right )^2$, the nature of the inequality in (\ref{eq:eq7}) will allow us to reduce the bound for $n$. Some experimenting shows that if we choose $C=6$, we can deduce that $n\leqslant 6$ using the same case checking procedure as before. \newline\newline Therefore, our strategy will be to find all the units of norm 1 which are between 1 and 6 (of which there are finitely many) and using (\ref{eq:eq1}) to look at the terms of the sequence up to the 30th and deduce the Zsigmondy bound. For $u>6$, we know from our above comments that the Zsigmondy bound is at most $6$, and there is little more we can say on this point. \newline\newline To find norm 1 units $1<u\leqslant 6$, we note that when $d\not \equiv 1 \pmod{4}$, $u$ is of the shape $u=a+b\sqrt{d}$, where $a,b$ are integers. Hence, the following inequality holds
\begin{displaymath}
2.618<a+b\sqrt{d}\leqslant 6\textrm{.}\end{displaymath}
Taking reciprocals we have
\begin{displaymath}
0.166<a-b\sqrt{d}<0.382 \textrm{,}\end{displaymath}
and it is clear that 
\begin{displaymath}
2\leqslant a \leqslant 3\textrm{.} \end{displaymath} 
If $a=2$ and $N(u)=1$ then we have $b^2d=3$. The only solutions of this are when $b^2=1$ and $d=3$ thus giving us $u=2\pm \sqrt{3}$. Hence, $u=2+\sqrt{3}$ is the only valid solution. Similarly if $a=3$ the only valid unit is $u=3+2\sqrt{2}$.\newline We now come to the case where $d\equiv 1 \pmod{4}$. A similar analysis for $u=\frac{a+b\sqrt{d}}{2}$ yields 
\begin{displaymath}
3\leqslant a \leqslant 6 \textrm{.}\end{displaymath}
The only solutions to $N(u)=1$ with $a$ in this range are $u=\frac{3\pm\sqrt{5}}{2}$ and $u=\frac{5\pm\sqrt{21}}{2}$, but again since $u>2.618$, we take the positive sign.
Hence there are four units of norm 1 which are greater than 1 but less than 6, namely $2+\sqrt{3}, 3+2\sqrt{2}, \frac{3+\sqrt{5}}{2}, \frac{5+\sqrt{21}}{2}$.
\newline\newline  We start with the case when $u=2+\sqrt{3}$, and we observe that for $7\leqslant n \leqslant 30$ inequality (\ref{eq:eq7}) holds when $n=8,10,12$. We also note that condition (\ref{eq:eq1}) fails when $n=8,10,12$ so $\Delta_8,\Delta_{10}$ and $\Delta_{12}$ all have PPDs, so we can restrict our attention to when $n\leqslant 6$. Here is a table illustrating the prime factors of $\Delta_n$ for $n$ from 1 to 6.
\begin{center}
\begin{tabular}[c]{|c|c|c|}
\hline
$n$ & $\Delta_n$ & Prime factors of $\Delta_n$ \\
\hline
1 & $-2$ & $2$ \\
\hline
2 & $-12$ & $2,3$ \\
\hline
3 & $-50$ & $2,5$ \\
\hline
4 & $-192$ & $2,3$ \\
\hline
5 & $-722$ & $2,19$ \\
\hline
6 & $-2700$ & $2,3,5$ \\
\hline
\end{tabular}
\end{center}
Therefore, the 4th and 6th terms of this sequence are the only ones which do not have a PPD.
\newline\newline We now turn our attention to $u=3+2\sqrt{2}$, where inequality (\ref{eq:eq7}) does not hold for any $n\geqslant 7$. So we can say immediately that $Z(\Delta)\leqslant 6$. We illustrate the prime factors of $\Delta _n$ in a table as previously.
\begin{center}
\begin{tabular}[c]{|c|c|c|}
\hline
$n$ & $\Delta_n$ & Prime factors of $\Delta_n$ \\
\hline
1 & $-4$ & $2$ \\
\hline
2 & $-32$ & $2$ \\
\hline
3 & $-196$ & $2,7$ \\
\hline
4 & $-1152$ & $2,3$ \\
\hline
5 & $-6724$ & $2,41$ \\
\hline
6 & $-39200$ & $2,5,7$ \\
\hline
\end{tabular}
\end{center}
It is therefore clear that when $u=3+2\sqrt{2}$ that all terms beyond the second have a PPD, so $Z(\Delta)=2$. 
\newline\newline When $u=\frac{5+\sqrt{21}}{2}$, we have that the inequality (\ref{eq:eq7}) holds when $n=12$ and that (\ref{eq:eq1}) is false when $n=12$. So we only need check terms from the sixth downwards to see which ones, if any, have no PPDs. We again list these terms and their prime factors in the table below.
\begin{center}
\begin{tabular}[c]{|c|c|c|}
\hline
$n$ & $\Delta_n$ & Prime factors of $\Delta_n$ \\
\hline
1 & $-3$ & $3$ \\
\hline
2 & $-21$ & $3,7$ \\
\hline
3 & $-108$ & $2,3$ \\
\hline
4 & $-525$ & $3,5,7$ \\
\hline
5 & $-2523$ & $3,29$ \\
\hline
6 & $-12096$ & $2,3,7$ \\
\hline
\end{tabular}
\end{center}
 Hence, we deduce again that $\Delta_n$ has a PPD for all terms beyond the sixth, and $\Delta_6$ is the only term which fails to have a PPD.
\newline\newline Finally, when $u=\frac{3+\sqrt{5}}{2}$ we find inequality (\ref{eq:eq7}) holds when $n=8,9,10,12,14,18,24,30$. Our condition (\ref{eq:eq1}) does not hold when $n$ is equal to 14, 18, 24, or 30. So we now need to check cases $n\leqslant 12$, to see which terms of $\Delta_n$ have primitive prime factors. These have all been listed in the table below.
 \begin{center}
\begin{tabular}[c]{|c|c|c|}
\hline
$n$ & $\Delta_n$ & Prime factors of $\Delta_n$ \\
\hline
1 & $-1$ & None \\
\hline
2 & $-5$ & $5$ \\
\hline
3 & $-16$ & $2$ \\
\hline
4 & $-45$ & $3,5$ \\
\hline
5 & $-121$ & $11$ \\
\hline
6 & $-320$ & $2,5$ \\
\hline
7 & $-841$ & $29$ \\
\hline 
8 & $-2205$ & $3,5,7$ \\
\hline 
9 & $-5776$ & $2,19$ \\
\hline
10 & $-15125$ & $5,11$ \\
\hline
11 & $-39601$ & $199$ \\
\hline
12 & $-103680$ & $2,3,5$ \\
\hline
\end{tabular}
\end{center}
It is immediately clear that $\Delta_n$ has no PPDs precisely when $n=6,10,12$.
\newline\newline We have therefore proven the following theorem
\begin{thm}
 Let $R$ be the ring of integers of the field $\mathbb{Q}(\sqrt{d})$, where $d$ is a squarefree positive integer. Let $0<u\in R$ be a unit of norm 1. If $u<\frac{1}{6}$ or $u>6$, then $Z(\Delta)\leqslant 6$. For all other such units $u$, one of the following holds:
\begin{itemize}
\item $u=3+2\sqrt{2},3-2\sqrt{2}$, where $Z(\Delta)=2$.
\item $u=2+\sqrt{3},2-\sqrt{3}$, where $Z(\Delta)=6$ and the only terms without a primitive prime divisor are $\Delta_4$ and $\Delta_6$.
\item $u=\frac{3+\sqrt{5}}{2}, \frac{3-\sqrt{5}}{2}$, where $Z(\Delta)=12$ and the only terms without a primitive prime divisor are $\Delta_6$, $\Delta_{10}$ and $\Delta_{12}$.
\item $u=\frac{5+\sqrt{21}}{2}, \frac{5-\sqrt{21}}{2}$, where $Z(\Delta)=6$ and the only term without a primitive prime divisor is $\Delta_6$.
\end{itemize}
\end{thm}

\subsection{Units of Norm $-1$}
We now wish to establish a similar result when $u$ is a unit of norm $-1$. If $u$ has norm $-1$ and $n=2k$ where $k$ is an odd integer, then (\ref{eq:eq1}) will not hold, but when $n$ is of this form $\Delta_n$ does not have a PPD, so we just ignore these values of $n$. Define $\Delta^{\prime}$ to be the sequence obtained by removing from $\Delta$, the terms $\Delta_n$ for which $n\equiv 2 \pmod{4}$. 
\newline\newline If $u$ is a quadratic unit of norm $-1$, then $u>\frac{1+\sqrt{5}}{2}$ and so by substituting into (\ref{eq:eq3}) we obtain that 
\begin{displaymath}
g(n)<3.44217\textrm{,} \end{displaymath}
where $g(x)$ is as before. Solving this inequality again using Maple 9.5 yields that \begin{displaymath}
n\leqslant 3375 \textrm{.} \end{displaymath}
Observing that if $u$ is norm $-1$, inequality (\ref{eq:eq8}) holds, so by Theorem 2.1 we have
\begin{equation} \label{eq:eq10}
u^{\phi(n)}<e^{5.03933}n^2\textrm{.} \end{equation}
Since $u\geqslant \frac{1+\sqrt{5}}{2}$, we are led to solve the following inequality 
\begin{displaymath}
\left ( \frac{1+\sqrt{5}}{2} \right ) ^{\phi(n)}<e^{5.03933}n^2\textrm{,} \end{displaymath}
checking cases on Maple 9.5 for $n$ between 7 and 3375  finds that this inequality is only true when 
\begin{displaymath}
n\leqslant 90\textrm{.} \end{displaymath}
Using the same trick as for norm 1, we observe that for $u>13$, inequality (\ref{eq:eq10}) implies that $n\leqslant 6$. We will therefore look at the cases $u>13$ and $u\leqslant 13$ separately. Finding the positive units of norm $-1$ that are between 1 and 13 is a finite problem and using the method from earlier we find that they are $1+\sqrt{2},\frac{1+\sqrt{5}}{2},2+\sqrt{5},\frac{11+5\sqrt{5}}{2},3+\sqrt{10},\frac{3+\sqrt{13}}{2},4+\sqrt{17},5+\sqrt{26},\frac{5+\sqrt{29}}{2},6+\sqrt{37},\frac{7+\sqrt{53}}{2},\frac{9+\sqrt{85}}{2}$.\newline\newline Now we need to look at the terms of the sequence for $1\leqslant n \leqslant 90$ where $n\not\equiv 2 \pmod{4}$ to see which have no PPDs. Doing the individual case checks as in the norm 1 case we find that when $u=1+\sqrt{2}$, inequality (\ref{eq:eq10}) holds when $n=7,8,9,10,11,12,14,15,16,18,20,21,22,24,26,28,30,36,42$. We are however ignoring the terms $\Delta_n$ for which $n\equiv 2 \pmod{4}$, so this leaves us to check the cases $n=7,8,9,11,12,15,16,20,21,24,28,36$. However (\ref{eq:eq1}) is violated, for all these values of $n$, therefore we conclude that $Z(\Delta^{\prime})\leqslant 4$. Once again we check to see if $\Delta_n$ has a primitive divisor for the relevant values of $n$ between 1 and 5. Again we illustrate the factors of $\Delta_n$ in tabular form
\begin{center}
\begin{tabular}[c]{|c|c|c|}
\hline
$n$ & $\Delta_n$ & Prime factors of $\Delta_n$ \\
\hline
1 & $-2$ & $2$ \\
\hline
3 & $-14$ & $2,7$ \\
\hline
4 & $-32$ & $2$ \\
\hline
5 & $-82$ & $2,41$ \\
\hline
\end{tabular}
\end{center} 
It is at once clear that $\Delta_4=\Delta_3^{\prime}$ is the only term of $\Delta^{\prime}$ without a PPD.
\newline\newline For $u=\frac{1+\sqrt{5}}{2}$, relations (\ref{eq:eq1}) and (\ref{eq:eq10}) are enough to ensure that for $n>6$, $\Delta_n$ has a PPD unless $n=12,20,24$. So we now need to check all the terms up to the 24th to see which ones have primitive prime divisors, and then we are done. Here is the table
 \begin{center}
\begin{tabular}[c]{|c|c|c|}
\hline
$n$ & $\Delta_n$ & Prime factors of $\Delta_n$ \\
\hline
1 & $-1$ & None \\
\hline
3 & $-4$ & $2$ \\
\hline
4 & $-5$ & $5$ \\
\hline
5 & $-11$ & $11$ \\
\hline
7 & $-29$ & $29$ \\
\hline 
8 & $-45$ & $3,5$ \\
\hline 
9 & $-76$ & $2,19$ \\
\hline
11 & $-199$ & $199$ \\
\hline
12 & $-320$ & $2,5$ \\
\hline
13 & $-521$ & $521$ \\
\hline 
15 & $-1364$ & $2,11,31$ \\
\hline
16 & $-2205$ & $3,5,7$ \\
\hline
17 & $-3571$ & $3571$ \\
\hline
19 & $-9349$ & $9349$ \\
\hline
20 & $-15125$ & $5,11$ \\
\hline
21 & $-24476$ & $2,19,211$ \\
\hline
23 & $-64079$ & $139,461$ \\
\hline
24 & $-103680$ & $2,3,5$ \\
\hline
\end{tabular}
\end{center}
We see that $\Delta_{12}=\Delta_9^{\prime},\Delta_{20}=\Delta_{15}^{\prime}$ and $\Delta_{24}=\Delta_{18}^{\prime}$ are the only terms of $\Delta^{\prime}$ that fail to have a PPD.
\newline\newline For all of the other units $u\leqslant 13$, conditions (\ref{eq:eq10}) and (\ref{eq:eq1}) are enough to secure that $n<6$, and hence that $Z(\Delta^{\prime})\leqslant 4$. Checking for primitive divisors of the remaining terms in exactly the same way as above, yields that all terms have a primitive prime divisor and so $Z(\Delta^{\prime})=1$.
\newline\newline Our case checking is now complete and we have derived the following result.
\begin{thm} Let $R$ be as in Theorem 3.1 and $1<u\in R$ be a unit of norm $-1$. Then for all $u>13$, $Z(\Delta^{\prime})\leqslant 4$. If $u\leqslant 13$, then one of the following is true:
\begin{itemize}
\item $u=1+\sqrt{2}$, where $Z(\Delta^{\prime})=3$ and the only term without a primitive prime divisor is $\Delta^{\prime}_3$;
\item $u=\frac{1+\sqrt{5}}{2}$, where $Z(\Delta^{\prime})=18$ and the only terms without a primitive prime divisor are $\Delta^{\prime}_{9}, \Delta^{\prime}_{15}$ and $\Delta^{\prime}_{18}$;
\item $Z(\Delta^{\prime})=1$.
\end{itemize}
\end{thm}

Combining the results of Theorems 3.1 and 3.2, gives us the statement of Theorem 1.4.

\end{document}